\documentclass[12pt]{article}

\usepackage{amssymb}
\usepackage{amsthm}
\usepackage{amsmath} 
\usepackage{stmaryrd}
\usepackage{mathrsfs,bm,hyperref}
\usepackage{soul,xcolor}

\usepackage{tikz}
\usetikzlibrary{shapes.geometric}
\usepackage{pgfplots}
\pgfplotsset{compat=1.18}
\usetikzlibrary{patterns}
\usetikzlibrary{perspective}
\usetikzlibrary{arrows}
\usetikzlibrary{3d,calc}
\usetikzlibrary{angles, calc, intersections, quotes, arrows.meta}
\usepackage{accents}
\usepackage{rotating}
\usepackage{soul}
\usepackage{ulem}
\usepackage{cite}

\usepackage{amsmath}
\usepackage{hyperref}
\usepackage{cleveref}
\usepackage{autonum}
\usepackage{mathtools}

\usepackage{dsfont}

\newcommand{\bk}{{\bm k}}

\renewcommand{\c}{\bm c}
\newcommand{\de}{\partial}
\newcommand{\dst}{\displaystyle}
\renewcommand{\emptyset}{\varnothing}
\newcommand{\ep}{\varepsilon}
\renewcommand{\leq}{\leqslant}
\renewcommand{\geq}{\geqslant}

\newcommand{\n}{\bm n}

\newcommand{\noin}{\noindent}

\newcommand{\om}{\omega}
\newcommand{\rf}[1]{(\ref{#1})}
\renewcommand{\rho}{\varrho}

\newcommand{\x}{{\bm x}}
\newcommand{\ym}{u_{-}}
\newcommand{\yp}{u_{+}}
\newcommand{\z}{{\bm 0}}
\newcommand{\up}{u^{+}}
\newcommand{\um}{u^{-}}

\DeclareMathAlphabet{\msfsl}{OT1}{cmss}{m}{sl}

\newcommand{\A}{\bm {\msfsl {A}}}
\newcommand{\B}{\bm {\msfsl {B}}}
\newcommand{\Bh}{\widehat{\B}}
\newcommand{\C}{{\mathbb C}}

\newcommand{\F}{\bm {\msfsl {\Phi}}}
\newcommand{\G}{\bm {\msfsl {G}}} 
\renewcommand{\H}{\bm {\msfsl {H}}}
\newcommand{\He}{\mathcal H_{\ep}}
\newcommand{\M}{\bm {\msfsl {M}}}
\newcommand{\N}{\bm {\msfsl {N}}}
\newcommand{\Np}{{\N}^{+}}
\newcommand{\Nm}{{\N}^{-}}

\newcommand{\Om}{\Omega}

\newcommand{\R}{{\mathbb R}}
\renewcommand{\Re}{\bm {\msfsl {R}}}

\newtheorem{theorem}{Theorem}[section]

\newtheorem{lemma}{Lemma}[section]

\newtheorem{remark}{Remark}

\newcommand{\I}{\ensuremath{\mathrm{i}}}
\newcommand{\E}{\ensuremath{\mathrm{e}}}
\newcommand{\D}{\ensuremath{\mathrm{d}}}

\begin{document}
\setstcolor{red}

\title{Spectral properties of operators and wave propagation \\ in high-contrast media\\ \vspace{5mm}}

\author{Yuri A. Godin$^a$\footnote{email: ygodin@charlotte.edu}, Leonid Koralov$^b$\footnote{email: koralov@umd.edu}, and Boris Vainberg$^a$\footnote{email: brvainbe@charlotte.edu} \vspace{5mm} \\
\small $^a$ Department of Mathematics and Statistics, \\ \small University of North Carolina at Charlotte,
Charlotte, NC, USA;  \\ \small $^b$ Department of Mathematics, University of Maryland, \\ \small College Park,
MD, USA}

\maketitle

\begin{abstract}
The paper aims to study the spectral properties of elliptic operators with highly inhomogeneous coefficients and related issues concerning wave propagation in high-contrast media.
A unified approach to solving problems in bounded domains with Dirichlet or Neumann boundary conditions, as well as in infinite periodic media, is proposed.
For a small parameter $\ep>0$ characterizing the contrast of the components of the medium, the analyticity of the eigenvalues and eigenfunctions is established in a neighborhood of $\ep =0$.
Effective operators corresponding to $\ep = 0$ are described.
\end{abstract}

\noin
{\bf Keywords} \\High-contrast medium; spectrum; analytic continuation;
 dispersion relation; \\ non-local boundary-value problem

\section{Introduction}
\label{sec1}

Problems involving media with high-contrast inclusions arise in many areas of physics and engineering: modeling flow in a porous medium
\cite{Chu:2010}, impedance tomography \cite{Borcea:96}, wave propagation in perforated elastic media \cite{Ammari:2009}, homogenization of composites \cite{Smyshlyaev:08, Cherdantsev:2012, Cherednichenko:2015}, the calculation of thermal fields in materials with highly conductive components \cite{Ewing:2009}, description of the bandgap structure of photonic and phononic crystals
\cite{Figotin:98, Figotin:2001, Lipton:2011, Lipton:17}.

Inhomogeneous materials are often used for controlling and manipulating wave propagation: preventing waves of specific frequencies from
propagating in certain directions, non-reciprocal wave propagation, etc. Such effects are more pronounced if the properties of the media (e.g., mass density, permeability, shear moduli, etc.) differ significantly in the matrix and inclusions. Thus, in such problems, a small parameter $\ep$ naturally arises, equal to the ratio of the material parameters of the constituent components. For example, $\ep$ could be the ratio of the diffusion coefficients or the square of the ratio of the wave propagation speeds in the host material and the inclusion.

The presence of a small parameter $\ep$ in the problem suggests using asymptotic methods; see \cite{Vishik:1960, Lions:1973, Papanicolaou:2011, Sanchez-Palencia:1980, Panasenko:80} for earlier works and \cite{Ammari:15, Cherednichenko:20, Weinstein:21, Gorb:21, Kiselev:22, Lim:24} for more recent results on the asymptotic analysis of solutions in high-contrast media. The operator $\A_\ep$ of the problem
is elliptic
in a bounded domain of $\mathbb{R}^d$ and has a discrete spectrum when $\ep > 0$. However, the domain of the operator
$\A_\ep$ depends on $\ep$ and the equation corresponding to $\A_\ep$ contains a factor $1/\ep$.
Thus, the operator $\A_\ep$ does not have a limit as $\ep \to 0$, and the asymptotic analysis of the problem is not quite obvious.

Despite this difficulty, we suggest a simple approach that allows us to prove the existence of the analytic extension of the eigenvalues
and eigenfunctions from $\ep >0$ to a neighborhood of $\ep =0$ and provide an explicit description of the eigenvalues
and eigenfunctions at $\ep = 0$.
We propose a unified approach that allows us to solve different high-contrast problems with the Dirichlet or Neumann boundary conditions, as well as periodic problems in $\mathbb{R}^d$ with the Bloch boundary conditions, where we study dispersion relations.

Different approaches were used
to derive the characteristic equation for the limit spectrum for the Dirichlet \cite{Panasenko:80} and Neumann \cite{Kiselev:22}
problems with an estimate of order $\ep$ between the spectrum of $\A_\ep$ and the limiting spectrum.
Some aspects of the analyticity of the Bloch spectrum were studied in \cite{Lipton:17}.

We employ two simple tricks to solve high-contrast problems with any boundary conditions and obtain the result using fairly simple arguments. First, we study the eigenvalues and eigenfunctions not for $\A_\ep$ but for the inverse operator $\B_\ep = \A^{-1}_\ep$.  Second, we study $\B_\ep$ by reducing the equation $\A_\ep u = f$ to an equation for the trace of $u$ on the boundary between the matrix and the inclusion. The latter equation is analytic in $\ep$, and therefore only general simple statements from functional analysis are needed.

The characteristic equation for the limit eigenvalues (as $\ep\downarrow 0$) has the form of an eigenvalue problem with Dirichlet data on the boundary of the inclusion equal to a constant determined by the equation containing the integral of the normal derivative of the solution
over the boundary of the inclusion. This type of non-local boundary-value problems have appeared in different publications; see, for example, \cite{Panasenko:80, Kiselev:22, FKW:2017}.

Diffusion processes leading to such non-local boundary value problems were considered in \cite{FKW:2017, FKW:2019, FK:2020}. In particular, in \cite{FKW:2017, FKW:2019} processes with non-local boundary behavior were shown to appear as limits of
processes with continuous trajectories (and, in the case of \cite{FKW:2019}, with diffusion coefficients that experience a sharp but continuous change near a hypersurface). Here, we obtain this non-local eigenvalue problem
in the case of a jump of the speed of the wave propagation (or the diffusion coefficient) as one of the results of our
general approach based on functional analysis.

The paper is organized as follows. For simplicity, we consider the Laplace operator with piecewise-constant coefficients
in sections \ref{D}-\ref{MI}. Sections \ref{D} and \ref{N} deal with the operator in a bounded domain with the Dirichlet and Neumann boundary conditions, respectively.
Bloch waves in periodic media in $\R^d$ are studied in section \ref{P}. The case of multiple inclusions is covered in section \ref{MI}. An extension of these results to elliptic operators with variable coefficients
 is presented in section \ref{VC}. Several simple examples are given in section \ref{E}.

\section{The Dirichlet problem}
\label{D}

Let $\Om\subset \R^d$ be a bounded domain containing a connected subdomain (inclusion) $\Om_-$. The case when $\Om_-$ has several components will be considered in section \ref{MI}.
We assume that the boundaries $\de\Om$ and $\Gamma= \de \Om_{-}$ are smooth enough (for example, $\de\Om,\Gamma \in C^{1,1}$, i.e., functions describing the boundaries have derivatives satisfying the Lipschitz condition),  and $\overline{\Om}_-\subset\Om$. We will write functions $u$ in $\Om$ in the form $u=(\up, \um)$, where $u^\pm$ are the restrictions of $u=u(\x)$ on the domains $\Om_\pm$, respectively, and $\Om_+=\Om\setminus\overline{\Om}_-$. Consider the Dirichlet problem for the negative Laplacian in $\Om$ with a piecewise-constant coefficient equal to $1, \ep^{-1},~0<\ep\ll 1,$ in $\Om_\pm$, respectively.

Let $\A_\ep$ be the following operator in $L_2(\Om)$:
 \begin{align}
 \label{AA}
\A_\ep u=(-\Delta \up, -\frac{1}{\ep} \Delta \um ), \quad u\in \He,
\end{align}
where the domain $\He = \He(\Om)$ of $\A_\ep$ consists of functions $u$ such that $u^\pm$ belongs to Sobolev spaces $H^2(\Om_\pm),$ respectively, and the following boundary conditions hold:
 \begin{align}
 \label{dirbc0}
 \up=0, \quad \x\in\de\Om; ~~ \quad \up=\um, \quad  \frac{\de \um}{\de \n} =\ep\, \frac{\de \up}{\de \n}, \quad \x \in  \Gamma,
\end{align}
where $\n$ is outward with respect to $\Om_{-}$, unit normal vector. The unit normal vector on $\de \Om$ also points outward from $\Om$.
Note that the boundary conditions on $\Gamma$ appear naturally when $\A_\ep$ is defined by the quadratic form involving the Laplacian with a piecewise-constant coefficient.

Two reasons make it difficult to use standard perturbation theory when studying the spectrum of $\A_\ep$ and solutions of the equation $\A_\ep u-\lambda u=f$ as $\ep\downarrow 0$: the domain of $\A_\ep$ depends on $\ep$, and the limit of $\A_\ep$ as $\ep\downarrow 0$ does not exist. Considering the inverse operator $\A_\ep^{-1}$ instead of $\A_\ep$ allows us to avoid these difficulties.

To study the equation
\begin{align}
 \label{eqA}
\A_\ep u=f, \quad  u\in \He(\Om), \quad  f=(f_+,f_-)\in L_2(\Om),
\end{align}
consider two separate Dirichlet problems:
\begin{align}
 \label{psi0}
 -\Delta \um &= \ep f_-, \quad \um \in H^2(\Om_-), \quad \left. \um \right|_{\Gamma} = \phi ,
\end{align}
\begin{align}
 \label{psiP0}
- \Delta \up  &=f_+, \quad \up \in H^2(\Om_+), \quad \left. \up \right|_{\de \Om} = 0 ,\quad \left. \up \right|_{\Gamma} = \phi,
\end{align}
with the same arbitrary function $\phi \in H^{3/2}(\Gamma)$.

Denote by  $\B_\ep$ the operator inverse to $\A_\ep$:
\[
\B_\ep=(\A_\ep)^{-1}: L_2(\Om) \rightarrow \He (\Om)
\]
Let us consider the range of the operator $\B_\ep$ in a wider $\ep$-independent space $H^2 (\Om_{\pm})$ that consists of functions $u=(\up,\um)$ such that $u^\pm\in H^2 (\Om_{\pm})$. Let $J_\ep: \He (\Om)\rightarrow  H^2 (\Om_{\pm})$ be the embedding operator, i.e., application of the operator $J_\ep$ to $ u\in \He (\Om)$ allows us to disregard conditions \rf{dirbc0}.
Denote
\begin{align}
 \widehat{\B}_\ep=J_\ep(\A_\ep)^{-1}: L_2(\Om) \rightarrow H^2(\Om_{\pm}).
\label{bhat}
\end{align}
Note that the spectrum of the operator $\A_\ep$ is discrete (as for any elliptic problem in a bounded region), and all its eigenvalues are positive since $\langle\A_\ep u,u\rangle$ is positive for all nontrivial $u\in \He(\Om).$

It is essential that $\A_\ep, \ep>0,$ and $\widehat{\B}_\ep, \ep>0,$ have the same eigenfunctions with the eigenvalues which are inverse to each other and positive since the eigenfunctions of operators $\B_\ep, \widehat{\B}_\ep, \ep>0,$ with positive eigenvalues coincide.

\begin{theorem}
\label{thm10}
The operator $\widehat{\B}_\ep$
exists for $0< \ep\ll1$ and can be extended analytically for $|\ep|\ll1$.  If
\begin{align}\label{B0f}
\widehat{\B}_0 f =u_0=(\up_0,\um_0) ,
\end{align}
then $\um_0=c_0$ and $\up_0$ is the solution of \rf{psiP0} with $\phi=c_0$, where constant $c_0$ is defined by the equation
\begin{align}\label{deqs0}
\int_\Gamma \frac{\de \up_0}{\de \n}\,\D S =\int_{\Om_-} f_- \,\D \x.
\end{align}
\end{theorem}
\begin{remark}
\label{r1}
Note that $\widehat{\B}_0 \neq J_0\A_0^{-1}$ since $\A_0$ does not exist.
\end{remark}
\begin{remark}
\label{r2}
Let us represent $\up_0$ as $c_0\up_1+\up_2$, where $\up_1$ is the solution of \rf{psiP0} with $\phi=1$ and $f_+=0$, and $\up_2$ is the solution of \rf{psiP0} with $\phi=0$. Then the left-hand side in \rf{deqs0} takes the form $c_0 a+b$, where by Green's identity
\begin{align}
 a=\int_\Gamma \frac{\de \up_1}{\de \n}\, \D S = -\int_{\Om_{+}} |\nabla \up_1 |^2 \, \D \x \neq0,
\end{align}
and
$b=\int_{ \Gamma} \frac{\de \up_2}{\de \n}\, \D S=\int_{\de \Om} \frac{\de \up_2}{\de \n}\, \D S + \int_{\Om_{+}} f_{+}\, \D \x.$
Thus $c_0$ is defined by \rf{deqs0} uniquely.
\end{remark}
To prove the theorem, we will need three simple lemmas. The first one is an observation, which does not need a proof.
\begin{lemma}
\label{L0}
Relations \rf{eqA} differ from Dirichlet problems \rf{psi0}, \rf{psiP0} with the same function $\phi$ on the boundary only by the presence in \rf{eqA} of the boundary condition for the derivatives of $(\up,\um)$ on $\Gamma$ described in the last condition in \rf{dirbc0}.
\end{lemma}

Let  $\N^\pm $ be the Dirichlet-to-Neumann (DtN) operators:
\begin{align}
 \N^\pm \phi = \left. \frac{\de u^\pm}{\de \n} \right|_{\Gamma},
\end{align}
where $u^\pm$ are the solutions of \rf{psi0}, \rf{psiP0} with $f_\pm=0$. We will use the notation $\M^-(\ep f_-),~\M^+f_+$ for the normal derivatives of $u^\pm$ when $\phi=0$ and $f_\pm\in L_2(\Om_\pm)$. Note that the operators
\begin{align}
\N^\pm: H^{3/2}(\Gamma)\to H^{1/2}(\Gamma), \quad \M^\pm: L_2(\Om_\pm)\to H^{1/2}(\Gamma)
\end{align}
are bounded, functions $\N^\pm\phi, \M^+f_+$ do not depend on $\ep$, and $\M^-(\ep f_-)=\ep\M^-( f_-)$ is proportional to $\ep$.
\begin{lemma}
\label{L1}
For each $f\in L_2(\Om)$ and $\ep>0$, the relation $u^\pm|_\Gamma=\phi$ is a one-to-one map between the solutions $u$ of
the problem \rf{eqA}
and the solutions $\phi$ of the equation
\begin{align} \label{Nphi}
(\N^--\ep \N^+)\phi=\ep \M^+f_+-\M^-(\ep f_-), ~\quad \phi\in H^{3/2}(\Gamma),
\end{align}
 extended to $u^\pm$ by \rf{psi0}, \rf{psiP0}.
\end{lemma}
\begin{proof}
The statement follows immediately from Lemma \ref{L0} since the last boundary condition in \rf{dirbc0} is equivalent to \rf{Nphi}.
\end{proof}

\begin{lemma} \label{L2}
The operator
$$\N=\N^--\ep \N^+: H^{3/2}(\Gamma)\to H^{1/2}(\Gamma)$$
is Fredholm when $|\ep|\ll 1$ and, as an operator on $L_2(\Gamma)$, is symmetric for all real $\ep$.
\end{lemma}

\begin{proof}
 The symmetry of operators $\N^\pm$ is an immediate consequence of Green's identities for problems \rf{psi0}, \rf{psiP0}. Hence,  $\N$ is symmetric. $\N^\pm$ are elliptic pseudo-differential operators of the first order (their symbols and ellipticity can be found in \cite{vg}). Since small perturbations preserve the ellipticity of $\Nm$, the operator $\N$ with small $|\ep|$ is elliptic and is therefore Fredholm.
\end{proof}
\begin{proof}[Proof of Theorem \ref{thm10}]
 We use Lemma \ref{L1} co construct $\B_\ep=\A_\ep^{-1},~\ep>0,$ and $\widehat{\B}_\ep = J_\ep \B_\ep$. Since problems \rf{psi0}, \rf{psiP0} are uniquely solvable,  \rf{psiP0}  does not depend on $\ep$, and solutions of \rf{psi0} are analytic in $\ep$, Lemma \ref{L1} implies that the theorem will be proved if we show that the operator $\F_\ep: L_2(\Om)\to H^{3/2}(\Gamma)$  mapping $f$ into the solution $\phi$ of \rf{Nphi} is analytic in $\ep, ~|\ep|\ll1 $, and $\F_0f=c_0$.

Let us represent the space $H^{3/2}(\Gamma)$ as a direct sum $S\oplus S^\perp$, where $S$ is the space of constant functions and $S^\perp$ is orthogonal in $L_2(\Gamma) $ subspace consisting of functions from $H^{3/2}(\Gamma)$ with zero integral. We write $\phi\in H^{3/2}(\Gamma)$ in vector form: $\phi=(\phi^c, \phi^\perp)$ where  $\phi^c=\int_\Gamma\phi\, \D S$ and $\int_\Gamma\phi^\perp\, \D S=0$, and then rewrite equation \rf{Nphi} as a system for $(\phi^c, \phi^\perp)$. Since the kernel and cokernel of the operator $\Nm$ are the space $S$, from Lemma \ref{L2} it follows that the system for $(\phi^c, \phi^\perp)$ has the following form
\begin{align}\label{syst}
 \left(
 \begin{array}{cc}
 0 & 0 \\[2mm]
  0 &  \Nm_{22}
 \end{array}
\right) \left(
 \begin{array}{cc}
 \phi^c\\[2mm]
  \phi^\perp
 \end{array}
\right)
-\ep\left(
 \begin{array}{cc}
  \Np_{11}  & \Np_{12}  \\[2mm]
  \Np_{21}  &  \Np_{22}
 \end{array}
\right) \left(
 \begin{array}{cc}
 \phi^c\\[2mm]
  \phi^\perp
 \end{array}
\right)= \ep \left(
 \begin{array}{cc}
 (\M f)^c\\[2mm]
  (\M f)^\perp
 \end{array}
\right),
\end{align}
where $ \Nm_{22} $ is the operator $\Nm$ restricted to $S^\perp$, and the right-hand side is the vector form of the right-hand side in \rf{Nphi}.

Since the operator $\Nm_{22}$ is invertible on $S^\perp$, the operator $(\Nm_{22}-\ep\Np_{22})^{-1}$ is analytic in $\ep$ for $|\ep|\ll1$. Thus, solving the second equation in \rf{syst} for $ \phi^\perp$ in terms of $\phi^c$ and $(\M f)^\perp$ we obtain
\begin{align} \label{phip}
\phi^\perp=\ep \G_\ep \phi^c+\ep \H_\ep (\M f)^\perp, ~\quad \G_\ep :S\to H^{3/2}(\Gamma), \quad \H_\ep: S^\perp\to H^{3/2}(\Gamma),
\end{align}
where operators $\G_\ep, \H_\ep$ are analytic in $\ep, ~|\ep|\ll1 $. Now, the first equation in \rf{syst} implies that
\begin{align} \label{phic}
(\Np_{11}+\ep\Np_{12}\G_\ep)\phi^c=-\ep\Np_{12} \H_\ep (\M f)^\perp-(\M f)^c.
\end{align}

The operator $\Np_{11}$ is the operator of multiplication by a constant $a\neq0$ defined in Remark~\ref{r2} after Theorem \ref{thm10}. Hence, the operator $(\Np_{11}+\ep\Np_{12}\G_\ep)^{-1}$ and, therefore, $\phi^c$ are analytic in $\ep,~|\ep|\ll1$. This and \rf{phip} imply the analyticity of $\F_\ep, |\ep|\ll1$. From \rf{phip} and \rf{phic} it also follows that $\F_0f=c_0$, where constant $c_0$ is defined by the equation
\begin{align} \label{f}
\Np_{11}c_0+(\M^+ f_+)^c+(\M^- f_-)^c=0.
\end{align}
The first term here is equal to $c_0\int_\Gamma\frac{\de \up_1}{\de \n}\, \D S$, the second coincides with $\int_\Gamma\frac{\de \up_2}{\de \n}\, \D S$, and the last one equals $\int_\Gamma\frac{\de \um}{\de \n}\, \D S=-\int_{\Om_{-}} f_{-}\, \D \x$. Thus, \rf{f} is equivalent to the equation \rf{deqs0} for $c_0$.

\end{proof}
Recall that the spectrum of the operator $\A_\ep$ is discrete, and all its eigenvalues are positive.
\begin{theorem}\label{t2}
 (i) The eigenvalues $\lambda=\lambda_j(\ep)$ and the eigenfunctions $u_{j}(\ep,\x)$ of the operator $\A_\ep,~\ep>0,$ can be enumerated in such a way that, for each $j$, $\lambda_j(\ep)\to\infty$ as $\ep\downarrow 0$ or $\lambda_j(\ep)$ and $u_{j}(\ep,\x)$ have analytic extensions in $\ep$ for $|\ep|\ll1$.

 (ii) For each compact set $\Lambda\subset\C$ that does not contain the limiting points $\lambda_j(0)$ of the eigenvalues, the resolvent
 \begin{align}
 \label{res}
 (\A_\ep - z)^{-1}:L_2(\Om)\to H^2 (\Om_\pm)
 \end{align}
 converges in the operator norm uniformly in $z\in\Lambda$ to the limit $\Re_z$, as $\ep\downarrow 0$, where
\begin{align}
 \label{res}
 \Re_z=\widehat{\B}_0({\bm {\msfsl {I}}}-z\widehat{\B}_0)^{-1}.
 \end{align}
\end{theorem}
\begin{proof}
The first statement follows immediately from Theorem \ref{thm10} and the
theorem \cite[Chap. VII, \S 3]{Kato} on the analyticity of the eigenvalues and
eigenfunctions of symmetric operators depending analytically on a parameter.
The latter theorem must be applied to $\widehat{\B}_\ep$. To justify the
second statement, one needs additionally to express the resolvent of $\A_\ep$ through the resolvent of the inverse operator:
\[
(\A_\ep - z)^{-1}=\widehat{\B}_\ep({\bm {\msfsl {I}}}-z\widehat{\B}_\ep)^{-1}, \quad ~z\neq \lambda_j(0),~~~|\ep|\ll1,
\]
and pass to the limit as $\ep\downarrow 0$.
\end{proof}
The next theorem provides the description of the limit set $\{\lambda_j, u_j\}$ of the eigenvalues and eigenfunctions of the operator $\A_\ep$ as $\ep\downarrow 0$.

\begin{theorem}\label{t3}
 A function $u$ is a limit eigenfunction of the operator $\A_\ep$, as $\ep\downarrow0$, with the limit eigenvalue $\lambda$, if

(i) $\um\equiv c_0,~\x\in\Om_-,$ with some constant $c_0$,

(ii) $\up $ is a solution of the problem
 \begin{align}
 \label{A0}
  -\Delta \up &=\lambda \up, \quad \up \in H^2(\Om_+), \quad \left. \up \right|_{\de \Om} = 0,\quad \left. \up \right|_{\Gamma} = c_0,
 \end{align}

(iii) the following relation holds
  \begin{align}\label{ex2}
\int_\Gamma \frac{\de \up}{\de \n}\,\D S+ c_0 \lambda |\Om_-|=0.
\end{align}
\end{theorem}
\begin{remark}
The integral in \rf{ex2} can be viewed as a function of  $c_0$ and $\lambda$ defined when a solution of \rf{A0} exists. If there are two solutions $u,~v$ of \rf{A0} with the same $c_0\neq0$, then $\lambda$ is an eigenvalue of the Dirichlet problem in $\Om_+$ with the eigenfunction $w=u - v$. From Green's formula applied to $u$ and $w$, it follows that
$\dst \int_\Gamma \frac{\de w}{\de \n}\,\D S =0$.
Thus, the validity of \rf{ex2} does not depend on the choice of solution of \rf{A0}.
\end{remark}
\begin{remark}
If $c_0=0$, then $\lambda$ is an eigenvalue of the Dirichlet problem \rf{A0} with an eigenfunction $u$ such that $\int_\Gamma \frac{\de u}{\de \n}\,\D S=0$. Thus, not an arbitrary eigenvalue of the Dirichlet problem \rf{A0} (with $c_0=0$) is a limit eigenvalue
 of $\A_\ep$. If $c_0\neq0$, then it can be canceled, i.e., one can put $c_0=1$ in (i)-(iii), and equation \rf{ex2} becomes a characteristic equation for limiting eigenvalues $\lambda$.
\end{remark}
\begin{proof}
Since the quadratic form $\langle \A_\ep u, u \rangle$ is bounded below for $0<\ep\ll1$, we have $\|\B_\ep\|\leq b<\infty$ when $0<\ep\ll1$. Thus, the limit eigenvalues of $\A_\ep$ are greater than or equal to $1/b$.
From Theorems \ref{thm10} and \ref{t2}, it follows that the limit $(\up_0,\um_0)$, as $\ep\downarrow0$, of an eigenfunction of $\A_\ep$ with a limit eigenvalue  $\lambda>0$ is an eigenfunction of the operator $\widehat{\B}_0$ with the
eigenvalue  $1/\lambda$. The converse is also true: each eigenfunction of the operator $\widehat{\B}_0$ is a limit eigenfunction of $\A_\ep$ with their eigenvalues inverse. Hence, the set of limit eigenvalues $\lambda$ and eigenfunctions $u_0=( \up_0,\um_0)$ of $\A_\ep$ satisfies \rf{B0f} with $f=\lambda( \up_0,\um_0)$,
where $( \up_0,\um_0)$ is described in Theorem \ref{thm10} and \rf{deqs0} is true. Since $\um_0 = c_0$ and $f_{-} = c_0 \lambda$, formula \rf{deqs0} is equivalent to \rf{ex2}.
\end{proof}

\section{The Neumann problem}
\label{N}

Let us consider operator \rf{AA} when the Dirichlet boundary condition on $\de\Om$ in \rf{dirbc0} is replaced by the Neumann condition
\[
\frac{\de \up}{\de \n}=0,\quad \x\in \de\Om.
\]
We will use the same notation for the operator $\A_\ep$ and its domain $\He(\Om)$ when the first condition in \rf{dirbc0} is the Neumann condition.

Since the operator $\A_\ep$ has a one-dimensional kernel and cokernel, both consisting of constant functions, the inverse operator $\B_\ep=(\A_\ep)^{-1},~\ep>0$, exists only if
\[
\A_\ep: \He^\perp (\Om) \to L_{2}^\perp(\Om)
\]
is considered as an operator between the spaces of functions from $\He$ and $L_2(\Om)$, respectively, orthogonal to constants. We will use the notation $\A_\ep^\perp$ for this operator to distinguish it from the similar operator without the orthogonality condition in the domain and the range.
We preserve the notation ${\B}_\ep$ for the inverse operator
\begin{align}
 \B_\ep = \left(\A_\ep^\perp \right)^{-1}:  L_{2}^\perp(\Om)\to \He^\perp (\Om).
\end{align}
We will consider the range of the operator $\B_\ep$ in a wider $\ep$-independent space $H^{2,\perp} (\Om_{\pm})$ that consists of functions $u=(\up,\um)\in H^2 (\Om_{\pm})$ that are orthogonal to constants. We denote
\begin{align}
 \widehat{\B}_\ep = J_\ep \B_\ep:  L_{2}^\perp(\Om)\to H^{2,\perp}(\Om_\pm),
 \label{BhatN}
\end{align}
where $J_\ep: \He^\perp (\Om) \rightarrow  H^{2,\perp} (\Om_{\pm})$ is the embedding operator, i.e., application of the operator $J_\ep$ to $ u\in \He^\perp (\Om)$ allows us to disregard conditions \rf{dirbc0}.

\begin{theorem}
\label{t4}
(i) In the case of the Neumann boundary condition on $\de\Om$, the inverse operator $\widehat{\B}_\ep$ exists for $0< \ep\ll1$ and can be extended analytically for $|\ep|\ll1$.

(ii) The operator $\widehat{\B}_0$ has the form
\[
\widehat{\B}_0 f =u_0=(\up_0,\um_0) ,
\]
where $\um_0=c_0,~\up_0$ is the solution of the problem
\begin{align}
 \label{pA}
- \Delta \up  &=f_+, \quad \up \in H^2(\Om_+), \quad \left. \frac{\de \up}{\de \n} \right|_{\de \Om} = 0 ,\quad \left. \up \right|_{\Gamma} = c_0,
\end{align}
and the constant $c_0$ satisfies $c_0=-|\Om_-|^{-1}\int_{\Om_+} u_0^+dx.$
\end{theorem}
\begin{remark}
Note that condition \rf{deqs0} holds automatically for solutions of \rf{pA} with arbitrary $c_0$ since the left-hand side in \rf{deqs0} for solutions of \rf{pA} equals $-\int_{\Om_+}f_+\D\x$ and we assume that $f\in L_{2}^\perp(\Om)$.
\end{remark}
\begin{proof}
The proof of Theorem \ref{t4} is somewhat similar to the proof of Theorem \ref{thm10}, with some steps being simpler. Lemma \ref{L0} remains valid. We slightly change the statement of Lemma \ref{L1}. We assume that $f\in L_{2}^\perp(\Om)$. Then the solution of equation $\A_\ep u=f,~ \ep>0,$ exists and is defined uniquely up to an arbitrary constant. The same is true for solutions of \rf{Nphi} when $|\ep|\ll1$. The lemma now states that \rf{Nphi} provides a one-to-one correspondence between the set of solutions of  $\A_\ep u=f,~ \ep>0,$ and the set of solutions of \rf{Nphi} with $\ep>0$ extended to $u^\pm$ by \rf{psi0}, \rf{psiP0} with the Neumann boundary condition on $\de\Om$ in \rf{psiP0}. However, the orthogonality of $\phi$ to constants does not imply that the particular solution of $\A_\ep u=f$ defined by \rf{psi0}, \rf{psiP0} with the Neumann boundary condition on $\de\Om$ in \rf{psiP0} is orthogonal to constants.

We look for a solution $\phi$ of \rf{Nphi}, with $f\in L_{2}^\perp(\Om)$, in the space $S^\perp\subset H^{3/2}(\Gamma)$ of functions orthogonal to constants. Then \rf{Nphi} is equivalent to the second equation in \rf{syst}
with $\phi^c = 0$.
Since the operator $\N^-=\N^-_{22}$ is invertible, it follows that $(\N^--\ep \N^+)^{-1}$ is analytic in $\ep$ when $|\ep|\ll1$, and the solution of \rf{Nphi} has the form $\phi=\Phi_\ep f$, where the operator $\Phi_\ep$ is analytic in $\ep, ~|\ep|\ll1,$ and $\phi=0$ when $\ep=0$. Thus, the particular solution $u=\widetilde{u}$ of the equation $\A_\ep u=f$ constructed by $\phi$ admits an analytic continuation in $\ep$ for $|\ep|\ll1$. Moreover, $\widetilde{u}= (\widetilde{u}^+_0, \widetilde{u}^-_0)$ when $\ep=0$, where $\widetilde{u}^-_0=0$,
while $\widetilde{u}^+_0$ satisfies \rf{pA} with $c_0 =0$.
All that remains is to add the constant $c_0$ to the solution $\widetilde{u}$ to make the average value of the solution $u_0 = (\up_0, \up_0)$ to be zero.
\end{proof}

\begin{theorem}\label{t5} Let $\A_\ep$ be the operator \rf{AA} with the Neumann boundary condition on $\de\Om$. Then

 (i) The eigenvalues $\lambda=\lambda_j(\ep)$ and the eigenfunctions $u_{j}(\ep,\x)$ of the operator $\A_\ep,~\ep>0,$ can be enumerated in such a way that $\lambda_j(\ep)\to\infty$ as $\ep\downarrow0$ or $\lambda_j(\ep)$ and $u_{j}(\ep,\x)$ have analytic extensions in $\ep$ for $|\ep|\ll1$.

 (ii) For every compact set $\Lambda\subset\C$ that does not contain the limit points $\lambda_j(0)$ of the eigenvalues (in particular, $0\notin \Lambda$), the resolvent
 \begin{align}
 \label{res6}
 (\A_\ep - z)^{-1}:L_2(\Om)\to H^2 (\Om_\pm)
 \end{align}
 converges in the operator norm uniformly in $z\in\Lambda$ to the limit $\Re_z$, as $\ep\downarrow 0$, where
\begin{align}
 \label{res61}
 \Re_z=\widehat{\B}_0({\bm {\msfsl {I}}}-z\widehat{\B}_0)^{-1} \quad  on \quad L_2^\perp(\Om), \quad and  \quad \Re_z c=-c/z.
 \end{align}
 (iii) The set of the limit eigenvalues and eigenfunctions $\{\lambda_j,u_j\}$ of $\A_\ep$, as $\ep\downarrow0$, with $\lambda_j>0$
 is defined as in Theorem \ref{t3}
 with the Dirichlet boundary condition on $\de\Om$ in \rf{A0} replaced by the Neumann condition.
\end{theorem}
\begin{remark}
The operators $\A_\ep$ and $\A_\ep^\perp$ have the same positive eigenvalues and the corresponding eigenfunctions, and the resolvent
$\Re_z,~z\neq 0,$ can be expressed through $\Re_z^\perp=(\A_\ep^\perp-z)^{-1}$ using simple shifts: $\Re_zf=\Re_z^\perp(f-c)+c/z,~z\neq 0,$ where
$c=\frac{1}{|\Om|}\int_\Om f\D \x$ is the average value of $f$.
Both statements of this remark are direct consequences of the solvability condition for the Neumann problem.
\end{remark}
\begin{proof}
Keeping in mind the remark above, the statements of the theorem follow from
Theorem \ref{t4} by using the same general arguments that were used to prove
Theorems \ref{t2} and \ref{t3}. One small additional step is needed to complete the proof of
the last statement of the theorem due to a small difference in the
definition of constants $c_0$ in Theorems \ref{thm10} and \ref{t4}. This difference implies that, following the same arguments, we will obtain the formulas \rf{A0}, \rf{ex2} for $u = \up$ with the Neumann boundary condition on $\de\Om$ in \rf{A0} and the integral over $\Gamma$ in \rf{ex2} replaced by $\int_{\Om_+} u^+\D \x$. However, Green's formula for the pair $u, v$, where $u=\up$ and $v=1$ implies that $\lambda\int_{\Om_+} u^+d\x=\int_\Gamma \frac{\de u^+}{\de \n}\,\D S.$

\end{proof}

\section{Waves in periodic media, dispersion surface}
\label{P}

\begin{figure}[th]
\begin{center}
\begin{tikzpicture}[scale=0.7,>=Stealth]
\def\l{1.5}
  \draw[fill=lightgray!10] (-\l,-\l,-\l) -- (-\l,-\l,\l) -- (\l,-\l,\l) -- (\l,-\l,-\l) -- cycle; 
  \draw[fill=lightgray!10] (-\l,-\l,\l) -- (-\l,-\l,-\l) -- (-\l,\l,-\l) -- (-\l,\l,\l) -- cycle; 
 \draw[fill=lightgray!10] (-\l,\l,-\l) -- (-\l,\l,\l) -- (\l,\l,\l) -- (\l,\l,-\l) -- cycle; 
   \draw[fill=lightgray!10] (-\l,-\l,\l) -- (\l,-\l,\l) -- (\l,\l,\l) -- (-\l,\l,\l) -- cycle; 
  \draw[fill=lightgray!10] (\l,-\l,\l) -- (\l,-\l,-\l) -- (\l,\l,-\l) -- (\l,\l,\l) -- cycle; 

\draw [->, thick] (0,0,0) -- (0,0,2.6*\l);
\node [above] at (0,0.2,2.6*\l) {\large  $x_1$};
\draw [->, thick] (0,0,0) -- (0,1.7*\l,0) node [right] {\large  $x_3$};
\draw [->, thick] (0,0,0) -- (1.7*\l,0,0) node [above] {\large  $x_2$};
 \def\eggheight{9mm}
 \begin{turn}{45}
  \path[ball color=white, opacity = 1]
  plot[domain=-pi:pi,samples=100,shift={(0,0)}]
  ({.8*\eggheight *cos(\x/4 r)*sin(\x r)},{-\eggheight*(cos(\x r))})
  -- cycle;
  \path[ball color=cyan!20, opacity = 0.5]
  plot[domain=-pi:pi,samples=100,shift={(0,0)}]
  ({.8*\eggheight *cos(\x/4 r)*sin(\x r)},{-\eggheight*(cos(\x r))})
  -- cycle;
  \end{turn}
  \node at (0,0) {\large $\Om_{-}$};
  \node [above] at (0.8,1.2,0) {\large$\Pi$};
\end{tikzpicture}
\end{center}
\caption{
The cell of periodicity $\Pi$ of a homogeneous medium with an inclusion $\Om_-$.
}
\label{cell}
\end{figure}
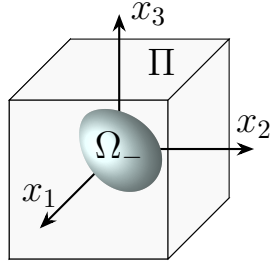

Consider wave propagation in $\mathbb{R}^d$ in the presence of periodically spaced inclusions, in which the wave speed is much higher than in the ambient medium. We assume that the fundamental cell of periodicity $\Pi$ is the cube $[-\pi,\pi]^d$. Let $\Om_-\subset \Pi$ denote the domain occupied by the inclusion in $\Pi$ (see Figure \ref{cell}).

 Waves in the whole space can be determined by the solutions $u=(\up,\um)$ of the following problem in $\Pi$:
\begin{align}
 \label{up}
 -\Delta \up= \lambda \up, \quad \x &\in \Pi \setminus \Om_-, \\[2mm]
 \label{um0}
 -\frac{1}{\ep}\,\Delta \um =  \lambda\um, \quad \x &\in \Om_-, \\[2mm]
 \label{cont0}
 \up = \um, \quad \frac{1}{\ep}\,\frac{\de \um}{\de \n} &= \frac{\de \up}{\de \n}, \quad \x \in \de \Om_-, \\[2mm]
 \label{B1}
 \left\rrbracket \up \E^{\I \bk \cdot \x} \right\llbracket &= 0.
\end{align}
Here $\bk\in \mathbb{R}^d$ is the wave (Bloch) vector, $\lambda=\om^2$, where $\om$ is the time frequency, $\n$ is the outward (with respect to $\Om_-$) unit normal vector, the coefficient $\ep^{-1} \gg 1$ in  \rf{um0} and \rf{cont0} reflects the contrast in the properties of the matrix and the inclusion, and the condition \rf{B1} means that $\up$ is sought in the Bloch form
\begin{align}
 \up(\x) = \Phi(\x) \E^{-\I \bk \cdot \x},
 \label{B}
\end{align}
where the function $\Phi$ is periodic with the period $2\pi$ in each variable, i.e., the inverted brackets $\rrbracket \cdot \llbracket $ denote
the jump of the enclosed quantity and its gradient across $\Pi$.

This problem differs from those considered in the previous sections only by the presence of the additional parameter $\bk$ in the boundary condition on $\de\Pi$. The eigenvalues $\lambda$ depend now on both $\ep$ and $\bk$. The multivalued function $\om=\om(\bk,\ep) = \sqrt{\lambda}$ for which problem \rf{up} - \rf{B1} has a nontrivial solution is called the {\it dispersion relation}.

We will focus on the generic case where the wave vector $\bk$ is not integer, that is, not all its components are integers. Then the Dirichlet-to-Neumann operator $\N^+=\N^+_{\bk}$ does not have a kernel, and problem \rf{up} - \rf{B1} is similar to the Dirichlet problem. If $\bk$ is an integer vector, problem \rf{up} - \rf{B1} can be easily studied similarly to the Neumann problem.

Denote by $\A_{\ep,\bk}, ~\ep>0,$ the operator \rf{AA} with $\Om=\Pi$. The space $\He(\Om)$ is now defined by \rf{dirbc0}, where the Dirichlet boundary condition on $\de\Om$ is replaced by the Bloch condition \rf{B1}. Denote $\B_{\ep,\bk}=\A_{\ep,\bk}^{-1}, ~\ep>0.$

We consider the range of the operator $\B_{\ep,\bk}$ in a wider $\ep$-independent space $H^2 (\Om_{\pm})$ that consists of functions $u=(\up,\um)\in H^2 (\Om_{\pm})$, $\Om_+ = \Pi \setminus \Om_{-}$. Let $J_\ep: \He (\Om)\rightarrow  H^2 (\Om_{\pm})$ be the embedding operator, i.e., application of the operator $J_\ep$ to $ u\in \He (\Om)$ allows us to disregard both conditions \rf{dirbc0} and \rf{B1}.
Denote
\begin{align}
 \Bh_{\ep,\bk}=J_\ep(\A_{\ep,\bk})^{-1}: L_2(\Om) \rightarrow H^2(\Om_{\pm}).
\end{align}

\begin{theorem}
\label{t6}
The operator $\Bh_{\ep,\bk}$ exists for $0< \ep\ll1$.  For any compact $K\subset \mathbb{R}^d$ that does not contain integer wave vectors $\bk$, there exists $\ep_0=\ep_0(K)>0$ such that the inverse operator $\Bh_{\ep,\bk}, ~ \ep>0,~\bk\in K,$ can be extended analytically in $\ep$ and $\bk$ for $|\ep|\leq \ep_0, ~\bk\in K$.  If
\begin{align}
\Bh_{0,\bk} f =(\up_\bk,\um_\bk) ,
\end{align}
then $\um_\bk=c_\bk,~\up_\bk$ is the solution of the problem
\begin{align}
 \label{upb}
- \Delta \up  &=f_+, \quad \up \in H^2(\Om_+), \quad  \left\rrbracket \up \E^{\I \bk \cdot \x} \right\llbracket = 0 ,\quad \left. \up \right|_{\Gamma} = c_\bk,
\end{align}
where $\Om_+=\Pi\setminus\Om_-,$ and the constant $c_\bk$ is determined from the equation
\begin{align}\label{d0}
\int_\Gamma \frac{\de \up_\bk}{\de \n}\,\D S =\int_{\Om_-} f_- \,\D \x.
\end{align}
\end{theorem}
\begin{proof}
The above statement is a complete analogue of Theorem \ref{thm10} and can be proved using the same arguments simply by replacing the Dirichlet boundary condition on $\de\Om$ in \rf{dirbc0}, \rf{psiP0} by \rf{B1}. The additional statement on analyticity in $\bk$ is also rather simple. Indeed, the solution of problem \rf{psiP0} with condition \rf{B1} on $\de\Om$ is analytic in $\bk\in\mathbb{R}^d$. The operators $\N^-, \M^-$ do not depend on $\bk$, and $\N^+, \M^+$ are analytic in $\bk\in\mathbb{R}^d$. The solution of system \rf{syst} is based on the invertibility of the operators  $\N^-_{22}, \N^+_{11}$ when $\bk\in K$. This leads to the analyticity of the solution $\phi$ of problem \rf{syst} and of $\Bh_{\ep,\bk}$.
\end{proof}
\begin{theorem}\label{t7}
 (i) For any compact $K\subset \mathbb{R}_\bk^d$ that does not contain integer wave vectors, the eigenvalues $\lambda=\lambda_j(\ep,\bk)$ and eigenfunctions $u_{j}(\ep,\x,\bk)$ of the operator $\A_{\ep,\bk},~\ep>0,~\bk\in K,$ can be enumerated in such a way that $\lambda_j(\ep,\bk)\to\infty$ as $\ep\to0$ or $\lambda_j(\ep,\bk)$ and $u_{j}(\ep,\x,\bk)$ have analytic extensions in $\ep$ for $|\ep|\ll1$.

 (ii) For any compact set $\Lambda \subset \C$ that does not contain the limit points $\lambda_j(0,\bk), \bk\in K,$ of the eigenvalues, the resolvent
 \begin{align}
 \label{res7}
 (\A_\ep - z)^{-1}:L_2(\Om)\to H^2(\Om_\pm),  \quad~z\in \C,
 \end{align}
 converges in the operator norm uniformly in $z\in\Lambda$ to a limit as $\ep\downarrow 0$.

 (iii) The set of limit eigenvalues  and eigenfunctions $\{\lambda_j,u_j\}$ of the operator $\A_{\ep,\bk},~\bk\in K,$ as $\ep\downarrow0$, is defined as in Theorem \ref{t3} with the Dirichlet boundary condition on $\de\Om$ in formula \rf{A0} replaced by \rf{B1}.
\end{theorem}
\begin{remark}
The last statement provides a description of the limit dispersion relation $\omega=\omega (\bk,0)$, $\bk\in K$, where $\omega = \sqrt{\lambda}$. A similar description for integer wave vectors $\bk$ can be obtained using the analogy with the Neumann problem.
\end{remark}
\begin{proof}
These statements follow from
Theorem \ref{t6} by repeating the same general arguments that were used to derive
Theorems \ref{t2} and \ref{t3} from Theorem \ref{thm10}.
\end{proof}

\section{Multiple inclusions}
\label{MI}

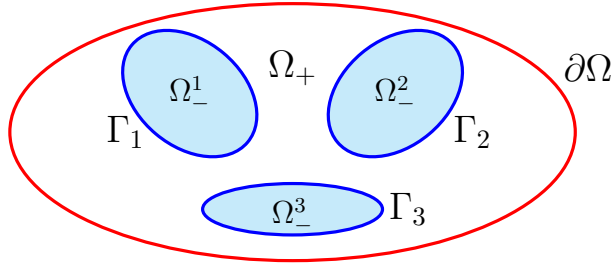
\begin{figure}[ht]
\begin{center}
\begin{tikzpicture}[scale=1.7,>=Stealth]

\draw[red, very thick] (0,0) ellipse (2.2cm and 1cm);

\node[draw, ellipse, blue, very thick, text = black, fill = cyan!20,
         minimum height=0.4cm, minimum width=2.0cm,rotate=40] (a) at (0.8,0.3) {{\rotatebox{-40}{$\Om_{-}^2$}}};

\node[draw, ellipse, blue, very thick, text = black, fill = cyan!20,
         minimum height=0.4cm, minimum width=2.0cm,rotate=-40] (b) at (-0.8,0.3)
         {{\rotatebox{40}{$\Om_{-}^1$}}};

\draw[red, very thick, blue,text = black, fill = cyan!20,] (0,-0.6) ellipse (0.7cm and 0.2cm);

\node at (0,0.5) {\mbox{\large $\Om_{+}$}};
\node at (2.3,0.5) {\large $\de \Om$};
\node at (0,-0.65) {$\Om_{-}^3$};
\node at (-1.3,0.0) {\large $\Gamma_1$};
\node at (1.4,0.0) {\large $\Gamma_2$};
\node at (0.9,-0.6) {\large $\Gamma_3$};

\end{tikzpicture}
\end{center}
\caption{
The domain $\Om$ in $\R^2$ containing three disjoint inclusions $\Om_{-}^1, \Om_{-}^2,$ and $\Om_{-}^3$ with
the boundaries $\Gamma_1, \Gamma_2,$ and $\Gamma_3$, respectively. $\Om_- =\cup_{i=1}^m \Om_{-}^i$,
$\Gamma =\cup_{i=1}^m \Gamma_i$, $m=3$.
}
\label{domain}
\end{figure}

We are going to discuss the changes needed to extend the previous results to the case when $\Om_- =\cup_{i=1}^m \Om_{-}^i$ with $\overline{\Om}_{-}^{\,i}$ disjoint and contained in $\Om$. We start with the Dirichlet problem discussed in Section \ref{D} (see figure \ref{domain}).

We represent the space $H^{3/2}(\Gamma)$ in the proof of Theorem \ref{thm10} as a direct sum $S\oplus S^\perp$, where $S$ is the space of functions that are constant on each $\Gamma_i$ and $S^\perp$ is orthogonal in $L_2(\Gamma) $ subspace consisting of functions from $H^{3/2}(\Gamma)$ with zero integrals over each $\Gamma_i$. We write $\phi\in H^{3/2}(\Gamma)$ in the vector form: $\phi=(\phi^c, \phi^\perp)$, where  $\phi^c=\frac{1}{|\Gamma_i|}\int_{\Gamma_i}\phi\, \D S,~\x\in\Gamma_i,$ and $\int_{\Gamma_i}\phi^\perp\, \D S=0, ~1\leq i\leq m$, and then rewrite equation \rf{Nphi} as a system \rf{syst} for $(\phi^c, \phi^\perp)$, where $\phi^c$ now is not a number but a vector $\c_0=(c_1,...,c_m)$. This leads to the version of Theorem \ref{thm10}, where $\um_0=c_i,~\x\in\Om_-^i,$ and $\up_0$ is the solution of \rf{psiP0} with $\phi=\phi^c=c_i,~\x\in\Gamma_i$. The following system of equations for the vector $\c_0$ now replaces equation \rf{deqs0}:
\begin{align}\label{deqs3}
\int_{\Gamma_i} \frac{\de \up_0}{\de \n}\,\D S =\int_{\Om_-^i} f_- \,\D \x, \quad 1 \leq i \leq m, \quad \up_0= \up_0 (\c_0,\x).
\end{align}

Let us show that \rf{deqs3} defines $\c_0$ uniquely. We represent $\up_0$ as
\begin{align}\label{sys}
 \up_0 = \sum_{i=i}^m c_i u_i (\x) + \tilde{u}(\x):=v(\x) + \tilde{u}(\x),
 \end{align}
where
 $u_i (\x)$ are the solutions of \rf{psiP0} with $f_+ = 0$ and $\phi = 1$ on
$\Gamma_i$ and $0$ for $\x \in \Gamma_j, ~j \neq i$, and $\tilde{u}(\x)$ is the solution of \rf{psiP0} with $\phi = 0$. This reduces \rf{deqs3} to a linear system for $c_i,~1\leq i\leq m$. If its determinant is zero, then there is a non-trivial vector $\c_0$ for which $\int_{\Gamma_i} \frac{\de v}{\de \n}\,\D S =0, ~ 1 \leq i \leq m.$ This and Green's formula for $v$ and $u_i$ imply that $\int_{\Om_+} \nabla v\cdot\nabla u_i \,\D \x =0, ~ 1 \leq i \leq m,$ and therefore $\nabla v = \z$
since $\nabla v= \sum_{i=1}^m c_i \nabla u_i$. Thus, $v$ is a constant on each connected component of $\Om_{+}$. If there are several components, they are separated by some domains $\Om^{i}_{-}$ (for example, $\Om$ is a ball $|\x|<1$
and $\Om_{-}$ is a spherical layer $1/3 < |\x| < 2/3$). All these constants are equal since
  $v = \up_0 - \tilde{u}$ is constant on each
$\Gamma_i = \de \Om_{-}^i$. (In the example above, $v$ is a constant on $\Gamma=\de\Om_-$, and therefore cannot have different constant values for $|\x|=1/3$ and $|\x|=2/3$.) Hence, $v$ is constant on $\Om_+$.
Thus, $v \equiv 0$ since $\left. v \right|_{\de \Om} =0$.
 This contradicts the fact that $v|_\Gamma=\c_0$ with a non-trivial $\c_0$. Hence, the determinant of the system \rf{deqs3} is not zero and \rf{deqs3} defines $\c_0$ uniquely.

Theorems \ref{t2}, \ref{t3} and their proofs remain valid with the description of the limit set $(\lambda_j, u_j)$ of the eigenvalues and eigenfunctions given by the following theorem.
 \begin{theorem}\label{tlast}
 A function $u=(\up,\um)$ is a limit eigenfunction of the operator $\A_\ep$, as $\ep\downarrow0$, with the limit eigenvalue $\lambda$, if

(i) $\um\equiv c_i,~\x\in\Om_-^i,~1\leq i\leq m,$ with some constant vector $\c_0=(c_1,...,c_m)$,

(ii) $\up $ is a solution of the problem
 \begin{align}
 \label{A05}
  -\Delta \up &=\lambda \up, \quad \up \in H^2(\Om_+), \quad \left. \up \right|_{\de \Om} = 0,\quad \left. \up \right|_{\Gamma_i} = c_i,\quad 1\leq i\leq m,
 \end{align}

(iii) the following relations hold
  \begin{align}\label{ex25}
\int_{\Gamma_i} \frac{\de \up}{\de \n}\,\D S+ c_i \lambda |\Om_-^i|=0,\quad 1\leq i\leq m.
\end{align}
\end{theorem}
Hence, if $u$ is a limit eigenfunction and $\c_0$ is the zero vector, then $\up$ is an eigenfunction of the Dirichlet problem \rf{A05} such that $\int_{\Gamma_i} \frac{\de \up}{\de \n}\,\D S=0,~1\leq i\leq m$.

Now, let $u$ be a limit eigenfunction with a limit eigenvalue $\lambda$ and a non-trivial vector $ \c_0$ in (i)-(iii). If $\lambda$ is not an eigenvalue of the Dirichlet problem \rf{A05} (with zero boundary conditions on all $\Gamma_i$), then the solution $\up$ of \rf{A05} depends linearly on $\c_0$, \rf{ex25} is a linear system for $\c_0$,  and equating the determinant of this system to zero yields a characteristic equation for the limit eigenvalues of $\A_\ep$ that are different from the eigenvalues of the Dirichlet problem \rf{A05}.

The results of this section are valid for the positive spectrum of the Neumann problem and for Bloch waves.

\section{Operators with variable coefficients}
\label{VC}

All the above results remain valid if the elliptic operator $\A_\ep$ has variable coefficients:
 \begin{align}
\A_\ep u=(-\nabla\cdot a^+(\x)\nabla \up, -\frac{1}{\ep}\nabla\cdot a^-(\x)\nabla \um ),  \quad a^\pm(\x)=(a^\pm_{ij}(\x)),
\quad 1 \leq i,j \leq d,
\end{align}
with the domain $\He(\Om)$ of $\A_\ep$ consisting of functions $u$ such that $u_\pm$ belongs to the Sobolev spaces $H^2(\Om_\pm),$ respectively, and the following boundary conditions hold:
 \begin{align}
 B\up=0, \quad \x\in\de\Om; ~~ \quad \up=\um, \quad  \frac{\de \um}{\de {\bm \nu}} =\ep\, \frac{\de \up}{\de {\bm \nu}}, \quad \x \in \Gamma,
\end{align}
where $B\up=0$ means the Dirichlet, Neumann (with the conormal derivative), or Bloch boundary condition, and $ \frac{\de u^\pm}{\de {\bm \nu}} = \langle  a^\pm(\x)  \nabla u,\n\rangle$ is the conormal derivative of $u^\pm$. One only needs to assume sufficient smoothness of the boundaries and coefficients to use the elliptic theory in the space $\He$. This smoothness assumption can be reduced by considering generalized solutions in the space $H^1_0$. The proofs do not need any changes except for using conormal derivatives instead of normal derivatives.

\section{Examples}
\label{E}
The examples below provide the limit spectrum for one-dimensional problems and a spherically symmetrical multidimensional problem. Statements concerning these examples can be obtained by applying Theorems
\ref{t3}, \ref{t5}, \ref{t7}, or by
solving the problems explicitly. These examples allow us to illustrate the results and, in some cases, to refine them.

Note that all the eigenvalues in the examples below are simple.

 \par{\it One-dimensional Dirichlet and Neumann problems.}
Consider the one-dimensional case on the interval $\Om = (-1,1)$ with the inclusion
$\Om_{-} = (a,b), ~ -1<a<b<1$. The Dirichlet eigenvalue problem for the operator
$\A_\ep$ (see \rf{AA},\rf{dirbc0}) has the form
\begin{align}
\label{sd}
-\yp^{\prime \prime} &=\lambda \yp, \quad x \in \Om_{+} =  (-1,a) \cup (b,1), \quad \left. \yp \right|_{|x|=1}=0,\\[2mm]
\label{sd1}
-\frac{1}{\ep}\ym^{\prime \prime} &= \lambda \ym, \quad x \in \Om_{-} = (a,b), \\[2mm]
\yp = \ym, \quad \ym^\prime &= \ep \yp^\prime, \quad x=a ~~\textrm{or}~~x=b.
\label{sd2}
\end{align}
The limit spectrum $S$ of this problem as $\ep \to 0$ consists of two parts, $S=S_1 \cup S_2$. If
$\frac{1+a}{1-b} = \frac{n_0}{m_0}$ is a rational number written as an irreducible fraction, then $S_1=\{ \lambda_{n,1}\}$, where
\begin{align}
 \lambda_{n,1} = \left(\frac{\pi m_0 n}{1-b}\right)^2 = \left(\frac{\pi n_0 n}{1+a}\right)^2, \quad n \geq 1.
\end{align}
The corresponding eigenfunction is proportional to
\begin{align}
 u_{n,1}(x) = \left\{
 \begin{array}{cc}
  (-1)^{n_0 n} \sin (\frac{\pi n_0 n}{1+a} (x+1)), & -1 \leq x \leq a, \\[2mm]
  0 & a \leq x \leq b, \\[2mm]
  (-1)^{m_0 n}\sin (\frac{\pi m_0 n}{1-b} (x-1)), & b \leq x \leq 1.
 \end{array}
\right.
\end{align}
If $\frac{1+a}{1-b}$ is irrational, then $S_1 = \emptyset$.

The elements of the set $S_2 = \{\lambda_{n,2}\}$ are the roots of the equation
\begin{align}
 \cot (\sqrt{\lambda_{n,2}}(1+a)) + \cot (\sqrt{\lambda_{n,2}}(1-b))  =\sqrt{\lambda_{n,2}} (b-a)
\end{align}
with the eigenfunctions proportional to
\begin{align}
\label{eigf1}
 u_{n,2}(x) = \left\{
 \begin{array}{cc}
  \frac{\dst \sin (\sqrt{\lambda_{n,2}} (x+1))}{\dst \sin (\sqrt{\lambda_{n,2}} (a+1))}, & -1 \leq x \leq a, \\[2mm]
  1 & a \leq x \leq b, \\[2mm]
  \frac{\dst \sin (\sqrt{\lambda_{n,2}} (x-1))}{\dst \sin (\sqrt{\lambda_{n,2}} (b-1))}, & b \leq x \leq 1.
 \end{array}
\right.
\end{align}

Next, consider a one-dimensional Neumann eigenvalue problem with the last condition in \rf{sd} replaced by $\left. \yp^\prime \right|_{|x|=1}=0$.

If
$\frac{1+a}{1-b}$ is a rational number that can be written in irreducible form as $ \frac{1+2n_0}{1+2m_0}$, then $S_1=\{ \lambda_{n,1}\}$, where
\begin{align}
 \lambda_{n,1} = \left(\frac{\pi (1+2nm_0)}{2(1-b)}\right)^2 = \left(\frac{\pi (1+2nn_0)}{2(1+a)}\right)^2, \quad n \geq 0.
\end{align}
The corresponding eigenfunction is proportional to
\begin{align}
 u_{n,1}(x) = \left\{
 \begin{array}{cc}
  (-1)^{n_0 n} \cos (\frac{\pi (1+2nn_0)}{2(1+a)} (x+1)), & -1 \leq x \leq a, \\[2mm]
  0, & a \leq x \leq b, \\[2mm]
  (-1)^{m_0 n+1}\cos (\frac{\pi (1+2nm_0)}{2(1-b)} (x-1)), & b \leq x \leq 1.
 \end{array}
\right.
\end{align}
If the condition on $\frac{1+a}{1-b}$ does not hold, then $S_1 = \emptyset$.

The elements of the set $S_2 = \{\lambda_{n,2}\}$ are the roots of the equation
\begin{align}
\label{tan}
 \tan (\sqrt{\lambda_{n,2}}(1+a)) + \tan (\sqrt{\lambda_{n,2}}(1-b))  +\sqrt{\lambda_{n,2}} (b-a)=0
\end{align}
with eigenfunctions defined by \rf{eigf1} with all the sine functions replaced by cosines and $\lambda_{n,2}$
 defined by \rf{tan}.
 \par{\it One-dimensional Bloch waves.}
 Here, the coordinates on the real line can be chosen so that $\Om_{-}=\{|x|<a<1\}$ and $\Pi=\{|x|<1\}$.
 Thus,
the one-dimensional version of problem \rf{up}-\rf{B1} has the form
\begin{align}
\label{y}
-\yp^{\prime \prime} &= \lambda \yp, \quad a < |x| < 1, \quad \yp=\Phi(x) \E^{-\I k x},\\[2mm]
-\frac{1}{\ep}\,\ym^{\prime \prime} &= \lambda\ym, \quad |x| < a, \\[2mm]
\yp = \ym, \quad \ym^\prime &= \ep \yp^\prime \quad |x| = a,
\label{ybc}
\end{align}
where $\Phi(x)$ is a periodic function with the period $p=2$ and $-\pi/2 < k \leq \pi/2$.
The eigenvalues are the roots of the equation
\begin{align}
 \cos (2\sqrt{\lambda_n}(1-a)) - \sqrt{\lambda_n} a \sin (2\sqrt{\lambda_n}(1-a)) = \cos 2k
\end{align}
with the eigenfunctions proportional to
\begin{equation}
 u(x) = \left\{
 \begin{array}{cc}
  \dst \frac{(1 + \E^{2\I k }) \cos (\sqrt{\lambda_n}(x+1))}{2\cos (\sqrt{\lambda_n}(1-a))} + \frac{(1 - \E^{2\I k })\sin (\sqrt{\lambda_n}(x+1))}{2\sin (\sqrt{\lambda_n}(1-a))}, & -1 \leq x \leq -a, \\[4mm]
  1 &  |x| \leq a, \\[2mm]
  \dst \frac{(1 + \E^{-2\I k })\cos (\sqrt{\lambda_n}(x-1))}{2\cos (\sqrt{\lambda_n}(1-a))}  - \frac{(1 - \E^{-2\I k })\sin (\sqrt{\lambda_n}(x-1))}{2\sin (\sqrt{\lambda_n}(1-a))} , \quad & a \leq x \leq 1.
 \end{array}
 \right.
\end{equation}

 \par{\it Concentric spheres in $\R^3$.}
Consider the Dirichlet problem
\begin{align}
 \label{up11}
 -\Delta \up= \lambda \up, \quad \x &\in \Om_{+} = \{a < |\x| < 1\}, \quad \left. \up \right|_{|\x|=1}=0,\\[2mm]
 \label{um011}
 -\frac{1}{\ep}\,\Delta \um =  \lambda\um, \quad \x &\in \Om_{-} = \{|\x|<a\}, \\[2mm]
 \label{cont011}
 \up = \um, \quad \frac{1}{\ep}\,\frac{\de \um}{\de r} &= \frac{\de \up}{\de r}, \quad |\x|=a.
\end{align}
The spherical symmetry and the connectivity of $\Om_{+}$ imply that the set $S_1$ is empty. The eigenvalues
from the set $S_2 = \{ \lambda_n\}$ are the roots of the equation
\begin{align}
 a\sqrt{\lambda_n} \cot (\sqrt{\lambda_n}(1-a)) = \frac{1}{3}\,\lambda_n a^2 -1
\end{align}
with eigenfunctions proportional to
\begin{align}
 u_n(r) = \left\{
\begin{array}{cc}
 1, & 0 \leq r \leq a, \\[2mm]
 \dst \frac{a \left(\sin (\sqrt{\lambda_n} r) \cos(\sqrt{\lambda_n} ) - \cos (\sqrt{\lambda_n} r) \sin(\sqrt{\lambda_n})\right)}
 {r \left(\sin (\sqrt{\lambda_n} a) \cos(\sqrt{\lambda_n} ) - \cos (\sqrt{\lambda_n} a) \sin(\sqrt{\lambda_n} )\right)}, & a \leq r \leq 1.
\end{array}
\right.
 \end{align}
Similar results are valid for the Neumann problem.

\section*{Funding}

The work of L. Koralov was supported by the NSF grant DMS-2307377 and the Simons Foundation grant MP-TSM-00002743.
The work of B. Vainberg was supported by the Simons Foundation grant 527180.

 \newcommand{\noop}[1]{}

\end{document}